\documentclass[11pt]{amsart}
\newtheorem{theorem}{Theorem}
\newtheorem{corollary}[theorem]{Corollary}
\newtheorem{example}[theorem]{Example}

\newtheorem{proposition}[theorem]{Proposition}

\begin{document}

\title{Normal conditional expectations of finite index \\ and sets of
modular generators}
\author{Michael Frank}
\address{Universit\"at Leipzig, Mathematisches Institut, D-04109 Leipzig, Germany}
\email{frank@mathematik.uni-leipzig.de}
\thanks{The author was partially supported by an N.S.F. grant}

\begin{abstract}
Normal conditional expectations $E: M \to N \subseteq M$ of finite index on
von Neumann algebras $M$ with discrete center are investigated to find an
estimate for the minimal number of generators of $M$ as a Hilbert $N$-module.
Analyzing the case of $M$ being finite type I with discrete center we obtain
that these von Neumann algebras $M$ are always finitely generated projective
$N$-modules with a minimal generator set consisting of at most $[K(E)]^2$
modular generators, where $[.]$ denotes the integer part of a real number and
$K(E) = \inf \{ K : \, K \cdot E - {\rm id}_M \geq 0 \}$. This result contrasts
with remarkable examples by P.~Jolissaint and S.~Popa showing the existence of normal
conditional expectations of finite index on certain type ${\rm II}_1$
von Neumann algebras with center $l_\infty$ which are not algebraically of
finite index, cf.~Y.~Watatani. We show that estimates of the
minimal number of modular generators by a function of $[K(E)]$ cannot exist
for certain type ${\rm II}_1$ von Neumann algebras with non-trivial center.
\end{abstract}
\maketitle



\smallskip \noindent
{\bf 1. Introduction. $\,$}
A conditional expectation $E: M \to N \subseteq M$ is said to
be of finite index if there exists a finite constant $K(E) = \inf \{ K : \,
K \cdot E - {\rm id}_M \geq 0 \}$, cf.~\cite{Popa95,FrKi98}. 
The goal of the present paper is the investigation of minimal generating
sets of right Hilbert W*-modules arising as $\{ M , E(\langle .,. \rangle_M) \}$
for conditional expectations $E: M \to N \subseteq M$ of finite index on von
Neumann algebras $M$ with discrete center $Z(M)$ (where $\langle x,y \rangle_M
= x^*y$ for every $x,y \in M$). We describe this case and indicate
partial solutions for the complementary case of von Neumann algebras with
diffuse center. The case of finite von Neumann algebras $M$ with discrete
center proves to be of major importance. Most interesting examples of normal
conditional expectations on von Neumann algebras of type ${\rm II}_1$ with
center $l_\infty$ that are non-algebraically of finite index were discovered
by P.~Jolissaint \cite[Prop.~3.8]{Jol91} and by S.~Popa \cite[Rem.~2.4]{Popa98}.
We give a detailed description of the situation and show the strong bounds
of this phenomenon. In particular, type ${\rm II}_1$ von Neumann algebras $M$
with non-trivial finite-dimensional center are always finitely generated
projective $N$-modules, however the minimal number of modular generators
cannot be estimated by a function of the characteristic constant $K(E)$, in
general.

\noindent
Analyzing the case of von Neumann algebras $M$ which are direct integrals of
a field of finite-dimensional type I factors over their discrete center we
obtain that these von Neumann algebras are always finitely generated projective
$N$-modules with a minimal modular generator set consisting of at most $[K(E)]^2$
generators, where $[.]$ denotes the integer part of a real number.

\medskip \noindent
Conditional expectations of finite index on von Neumann algebras have
remarkable properties:
they are automatically normal and faithful (\cite[Prop.~1.1]{Popa98}), the von
Neumann algebra $M$ is complete with respect to the Hilbert norm derived from
the new $N$-valued inner product $E(\langle .,. \rangle_M)$
(\cite[Prop.~3.3]{BDH88}), and the center of $M$ is finite-dimensional if and
only if the center of $N$ is. Moreover, they commute with the abstract
projections of von Neumann algebras to most of their type components, especially
to the parts with discrete and diffuse center (\cite[Th.~1]{FiIs97}) and to their
type ${\rm I}_{fin}$, ${\rm I}_\infty$, ${\rm II}_1$, ${\rm II}_\infty$ and
${\rm III}$ parts (\cite[1.1.2, (iii)]{Popa95}, \cite[\S 2]{FrKi98}). This
justifies the canonical character of our basic question for minimal generating
sets of $M$ as an $N$-module and opens up the chance to solve it in the form
of separate case studies of the different types of von Neumann algebras $M$
that can occur.

\smallskip \noindent
The case of properly infinite von Neumann algebras $M$ was solved by several
authors without any reference to special structures of the center:
there exists a single element $m \in M$ such that every element $x \in M$
can be represented as $x = mE(m^*x)$ and $E(m^*m) = 1_M$,
cf.~\cite[Lemma 3.21, Prop.~3.22]{BDH88}, \cite[Th.~1.1.6, (b)]{Popa95}.
The proofs show the existence of $m \in M$. We provide a constructive proof
of the generating element $m$ for properly infinite discrete type ${\rm I}$
von Neumann algebras during the proof of Theorem \ref{th-findiscr} below,
cf.~the final remarks.

\bigskip \noindent
{\bf 2. The type ${\rm II}_1$ case. $\,$}
Due to the type splitting properties of $E$ we can concentrate our efforts
on the finite case in the sequel. If $M$ and $N$ are type ${\rm II}_1$ factors
then $M$ is a finitely generated projective $N$-module with a minimal set of
modular generators consisting of at most $[K(E)]+1$ elements,
cf.~\cite[Prop.~1.3]{PiPo}. For type ${\rm II}_1$ von Neumann algebras
with discrete center finite generation of $M$ as an $N$-module can only be
guaranteed over finite-dimensional pieces of the center of $N$,
cf.~\cite{ACS}.

\begin{proposition}  \label{prop1}
  If $E: M \to N \subseteq M$ is a conditional expectation of finite
  index on a type ${\rm II}_1$ von Neumann algebra $M$ with discrete
  center then for any finite number of minimal central projections of $N$ with
  sum $p \in Z(N)$ the right Hilbert $pN$-module $Mp$ is finitely generated.
\end{proposition}

\begin{proof}
By \cite{BDH88,FrKi98} the von Neumann algebra $N=E(M)$ and
the set ${\rm End}_N(M)$ of all bounded $N$-linear operators on the Hilbert
$N$-module $\{ M , E(\langle .,. \rangle_M) \}$, which is a von Neumann
algebra, have to be of type ${\rm II}_1$ with discrete center. Both these
W*-algebras share a common center $Z(N)$ if we identify appropriate multiples
of the identity operator on $M$ with the corresponding elements of $Z(N)$.
Selecting a finite sum $p$ of minimal central projections of $Z(N)$ the
corresponding von Neumann algebra $p \cdot {\rm End}_N(M)$ does not possess
any non-unital two-sided norm-closed ideal. In particular,
\[
  p \cdot {\rm End}_N(M) \equiv {\rm lin} \{  p \cdot \theta_{x,y} \, : \,
  \theta_{x,y}(z)= y E(x^*z) \quad {\rm for} \quad z \in M \} \, .
\]
By \cite[Ex.~15.O]{NEWO} the Hilbert $pN$-module $\{ M, E(\langle .,.
\rangle_M) \}$ has to be finitely generated.
\end{proof}

\medskip \noindent
Surprisingly, the statement on type ${\rm II}_1$ von Neumann algebras $M$ to
be finitely generated projective $N$-modules for conditional expectations
$E: M \to N \subseteq M$ of finite index, in general fails to be
true, (cf.~\cite[Prop.~1]{FiIs97}). A remarkable counterexample to this 
conjecture was brought to my attention by S.~Popa, cf.~\cite[Rem.~2.4]{Popa98}.
A similar construction with a different proof can be found in a paper by
P.~Jolissaint \cite[Th.~2.2, Prop.~3.8]{Jol91}.
The counterexample also shows in its generalized version that an analogue of
Proposition \ref{prop1} cannot be true for type ${\rm II}_1$ von Neumann
algebras with diffuse center even if it is generated by a non-discrete finite
measure space.

\begin{example}      {\rm
Let $R$ be the hyperfinite type ${\rm II}_1$ W*-factor and select a projection
$p_k$ with ${\rm tr}(p_k)=2^{-k}$. Since the type ${\rm II}_1$ factors
$p_kRp_k$ and $(1-p_k)R(1-p_k)$ are $*$-isomorphic we can find a $*$-isomorphism
$\theta_k$ of them. Set $M_k = \{ x+y : x \in p_kRp_k , \, y \in (1-p_k)R(1-p_k)
\}$ and $N_k = \{ x + \theta_k(x) : x \in p_kRp_k \}$. Note that $p_kM_kp_k =
N_kp_k$, $(1-p_k)M_k(1-p_k) = N_k(1-p_k)$ and $N_k' \cap M_k = {\bf C}p_k \oplus
{\bf C}(1-p_k)$.


\noindent
Considering the unique trace-preserving conditional expectation
$E_{k,tr}: M_k \to N_k \subset M_k$ with $k \in {\bf N}$ we obtain
${\rm Ind}(E_{k,tr})= {\rm tr}(p_k)^{-1} + {\rm tr}(1-p_k)^{-1} =
2^k + 2^k/(1-2^k)$ by \cite[\S 2.5, Th., (iii)]{Jones88},
\cite[Lemma 2.2.2]{Jones83}. Following \cite[Prop.~3.15]{BDH88} any other
conditional expectation $E_k: M_k \to N_k \subset M_k$ can be expressed
as $E_k(x) = E_{k,tr}(a^*xa)$ for some element $a \in N_k' \cap M_k$.
Among these conditional expectations we can find one of minimal index,
$E_{k,min}$, with index value ${\rm Ind}(E_{k,min})=4$ independent of $k$.
Two specific values are $E_{k,min}(p_k)=E_{k,min}(1-p_k)=1/2$.

\noindent
Obviously, $M_k$ is a finitely generated projective $N_k$-module for every
$k \in {\bf N}$. However, the number of modular generators has to be greater
than $2^k$ since every generator $u_k \in M_k$ can contribute at most the
value $1$ to the index value ${\rm Ind}(E_{k,tr})= 2^k + 2^k/(1-2^k)$ of the
trace-preserving conditional expectation $E_{k,tr}$ because of the finiteness
of $M_k$, cf.~\cite[Th.~3.5]{BDH88}.
Indeed, ${\rm tr}(u_ku_k^*)={\rm tr}(u_k^*u_k)={\rm tr}(E_{k,tr}(u_k^*u_k))
\leq 1$ since the generators $u_k$ can be selected such that $E_{k,tr}
(u_k^*u_k)$ is a projection by \cite{BDH88}, and so the spectrum of $u_ku_k^*$
is restricted to a subset of the unit interval $[0,1]$, i.e.~$u_ku_k^* \leq 1$.

\noindent
In a final step we form W*-algebras $M$ and $N$ with a common direct integral
decomposition over $l_\infty({\bf N})$ attaching to every
minimal projection $q_k$ of $l_\infty({\bf N})$ either $M_k$ or $N_k$,
respectively.
The unique trace-preserving conditional expectation $E_{tr}: M \to N \subset M$
is not of finite index any longer, whereas the minimal conditional expectation
$E_{min}: M \to N \subset M$ has the index value ${\rm Ind}(E_{min})=4$.
Since the number of generators of $M_k$ as a finitely generated projective
$N_k$-module is greater than $2^k$ the total number of generators of $M$ as
a self-dual Hilbert $N$-module is infinite.

\noindent
We can find analogous examples for type ${\rm II}_1$ von Neumann algebras $M$
with diffuse center. To see this take the diffuse von Neumann algebra
${\rm L}^\infty([0,1],\lambda)$ and form two direct integrals over $[0,1]$
in such a way that all fibres on the intervals $[1/(k+1),1/k]$ are $*$-isomorphic
to $M_k$ or $N_k$, respectively, $(k=1,2,...)$. Denote the resulting von Neumann
algebras by $M$ and $N$, respectively. The von Neumann algebra $M$ cannot be
a finitely generated Hilbert $N$-module for the same reasons as for the
example with discrete center.      }
\end{example}

\noindent
The index value $4$ is not essential and can be replaced by any other number
within the range $[4,\infty)$ since we can select another conditional
expectation $E$ with a fixed finite index between ${\rm Ind}(E_{tr})=(2^{2k}-
2^{k+1})/2^k-1 \stackrel{k \to \infty}{\longrightarrow} \infty$ and ${\rm Ind}
(E_{min})=4$ relying on suitable elements in the non-trivial relative
commutants. The number of modular generators is not affected by this choice.
For index values lower than four the relative commutant has to be trivial, and
we are again in the factor case.

\begin{corollary}
Let $E: M \to N \subseteq M$ be a (normal) conditional expectation of finite
index on a W*-algebra $M$ with image algebra $N$. Then $M$ is not necessarily
a finitely generated projective $N$-module, especially if $N$ is of type
${\rm II}_1$ with infinite-dimensional discrete or arbitrary diffuse center.
\newline
Furthermore, even if $M$ has a finite-dimensional non-trivial center the
minimal number of ge\-nerators of $M$ as a finitely generated projective
$N$-module cannot be estimated by a fixed function of $K(E)$, in general.
\end{corollary}

\noindent
{\bf 3. The discrete finite case. $\,$}
The remaining case is that of a von Neumann algebra $M$ of type ${\rm I}$
with discrete center, but without properly infinite part in its central
direct integral decomposition. The situation for finite-dimensional centers
was partially considered by T.~Teruya \cite{Ter92}. For our purposes the
center of $M$ can be arbitrarily large. Also, the relative position of the
centers $Z(M)$ and $Z(N)$ can be arbitrary, i.e.~somewhere between the
simplest case $Z(M)=Z(N)$ and the most complicated case $Z(M) \cap Z(N) =
{\bf C}1_M$, see \cite{FrKi98} for examples.

\noindent
\begin{theorem} \label{th-findiscr}
Let $M$ be a discrete W*-algebra without properly infinite part. Consider a
(nor\-mal) conditional expectation $E: M \to N \subseteq M$ of finite
index possessing the structural constant $1 \leq K(E) = \inf \{ K : \, K
\cdot E - {\rm id}_M \geq 0 \} \in {\bf R}$.
Then $M$ is a finitely generated projective $N$-module with a minimal modular
generator set consisting of at most $[K(E)]^2$ generators, where $[.]$ denotes
the integer part of a real number. \newline
If $n=K(E)$ is an integer then the normalized trace $E$ on the full $n \times
n$-matrix algebra
\linebreak[4]
$M={\rm M}_n({\bf C})$ realizes the indicated upper bound $n^2$ for the number
of modular generators, i.e.~the estimate can in general not be improved.
\end{theorem}

\begin{proof}
We note in passing that the construction principles do not depend on whether
we consider $M$ as a right or left $N$-module. To be precise we
use the right $N$-module notation subsequently. Moreover, the considerations
below work equally well for the properly infinite discrete case for which a
stronger result is already known. So they are only of interest as a way to
construct the generating element $m \in M$ for this case.

\noindent
The strategy of the presented proof is to find a suitable set of modular
generators in the beginning, and to reduce this set of generators until we
obtain a sufficiently small generating set fulfilling the prediction.  Note
that one or the other reduction step might be redundant when we resort to
particular examples, however the stressed for statement requires a proof from
a general viewpoint.

\smallskip \noindent
Consider a minimal (with respect to $N$) projection $p \in N$. Let $q \leq p$
be a projection of $M$. Since $E$ is faithful $E(q) \not= 0$ and the
inequality $0 < E(q) \leq E(p) = p$ holds. Because of the minimality of $p$
inside $N$ there exists a number $\mu \in (0,1]$ such that $E(q)= \mu p$.
That is,
  \[
  E(\mu^{-1} \cdot q)=p \, .
  \]
Moreover, since $E$ was supposed to be of finite index with structural
constant $K(E) \geq 1$ we obtain
  \[
  K(E) \cdot p = K(E) \cdot E(\mu^{-1} q) \geq \mu^{-1} q > 0
  \]
and hence, the estimate $\mu \geq K(E)^{-1}$.
Suppose, $p \in N$ can be decomposed into a sum of pairwise orthogonal
(arbitrary) projections $\{ q_\alpha : \alpha \in I \}$ inside $M$.
Obviously, $q_\alpha \leq p$ for every $\alpha \in I$, and
  \[
  p = E(p) = \sum_{\alpha \in I} E(q_\alpha)
  = \left( \sum_{\alpha \in I} \mu_\alpha \right) p
  \geq \left( \sum_{\alpha \in I} K(E)^{-1}_{(\alpha)} \right) p \, .
  \]
Consequently, the sum has to be finite and the maximal number of non-trivial
summands is $[K(E)]$, the integer part of $K(E)$. We see, for the minimal
projections $p \in N \subseteq M$ every family of pairwise orthogonal
subprojections $\{ q_i \} \in M$ of $p$ is finite.

\smallskip \noindent
In this way we obtain that every minimal projection $p \in N$ has to be
represented only in a small part of the central direct integral decomposition
of $M$ which has a finite-dimensional center. Extending a minimal projection $p \in
N \subseteq M$ by the partial isometries of $N \subseteq M$ every W*-factor
block of the central direct integral decomposition of $N$ turns out to be
represented in a part of the central direct integral decomposition of $M$
with a finite-dimensional center. Conversely,
the conditional expectation of finite index $E$ maps each W*-factor block of
the central direct integral decomposition of $M$ into a part of the central
direct integral decomposition of $N$ with a finite-dimensional center. To see
that consider a minimal central projection $q \in M$ together with the
conditional expectation $E':qM \to qN \subseteq qM$, where $E'(x)=E(x) \cdot
E(q)^{-1} \cdot q$ for $x \in qM$. Since $K(E) \cdot E(q) \geq q$ by
supposition we obtain
\[
K(E) \cdot E(x) \cdot E(q)^{-1} \cdot q \geq K(E) \cdot E(x) \cdot q
        \geq qxq = xq
\]
for any $x \in qM$. Therefore, $K(E') \leq K(E)$, and ${\rm dim}(Z(qN)) \leq
[K(E)]$ by \cite[Th.~3.5, (ii)]{BDH88} and by the structure of a $qN$-module
basis of the Hilbert $qN$-module $\{ qM, E'(\langle .,. \rangle) \}$.

\medskip \noindent
To find a suitable set of generators of $M$ as a right Hilbert $N$-module
we decompose the identity $1_N=1_M$ into a w*-sum of pairwise orthogonal
minimal projections $\{ p_\nu \} \subset N$, and further into a suitable
subdecomposition of this sum into a w*-sum of pairwise orthogonal minimal
projections $\{ q_\alpha : \alpha \in I \} \subset M$. Without loss of
generality we will always assume that minimal projections $p_{\nu_1},p_{\nu_2}
\in N$ of our choice with $p_{\nu_1}=v^*v$, $p_{\nu_2}=vv^*$ for $v \in N$
possess finite sum decompositions $p_{\nu_1}=\sum_i q_{i,\nu_1}$,
$p_{\nu_2}=\sum_i q_{i,\nu_2}$ with finite sets of minimal projections
$\{ q_{i,\nu_1} \}$, $\{ q_{i,\nu_2} \} \subset M$ such that every partial
isometry $q_{i,\nu_1}v \in M$ has the corresponding domain projection
$q_{i,\nu_2}$.

\noindent
In our special setting a suitable $N$-module basis of $M$ contains this
maximal set of pairwise orthogonal minimal
projections $\{ q_\alpha \}$ of $M$ scaled down by the inverse of the
number $\mu_\alpha$ arising from the equality $E(q_\alpha) = \mu_\alpha p_\nu$
for some minimal projection $p_\nu \in N$ of our initial choice and certain
$\mu_\alpha \in [K(E)^{-1},1]$. As an intermediate result we get
\[
\{ \mu_\alpha^{-1/2} \cdot q_\alpha : \alpha \in I \} \subseteq
{\rm basis} \, .
\]
If $M$ is commutative, then the Hilbert $N$-module basis of $M$ is complete,
and
\[
{\rm Ind}(E) = \sum_{\alpha \in I} \mu_\alpha^{-1/2} \cdot q_\alpha
                     \cdot (\mu_\alpha^{-1/2} \cdot q_\alpha)^*
            \leq K(E) \sum_{\alpha \in I} q_\alpha
            \leq K(E) \cdot 1_M
\]
by \cite[Thm.~3.5]{BDH88}. However, if $M$ is non-commutative, then we have to
add all those minimal partial isometries $\{ u_\beta : \beta \in J \}$ of $M$
that connect two minimal projections of our choice $\{ q_\alpha \}$,
but also scaled down by
the inverse of the number $\mu_\beta$ arising from the equality
$E(u_\beta^*u_\beta) = \mu_\beta p_\nu$ for some minimal projection $p_\nu
\in N$ of our initial choice and certain $\mu_\beta \in [K(E)^{-1},1]$,
cf.~\cite[Ex.~1.1]{FrKi98}. Finally,
  \begin{equation} \label{mark1}
     \{ \mu_\alpha^{-1/2} \cdot q_\alpha : \alpha \in I \} \cup
     \{ \mu_\beta^{-1/2} \cdot u_\beta  : \beta \in J \}
     \equiv {\rm basis} \, .
  \end{equation}
Finally a Hilbert $N$-module basis of $M$ is complete, but often rather large.

\smallskip \noindent
In a next step we use Hilbert $N$-module isomorphisms to reduce the number
of generators in the generating set (\ref{mark1}). For simplicity we consider
the projections $\{ q_\alpha \}$ as special partial isometries, too.
Define equivalence classes of partial isometries of (\ref{mark1}) by the rule:
$\, u_\beta \sim u_\gamma$ if and only if $q_\alpha u_\beta \not= 0$,
$q_\alpha u_\gamma \not= 0$ for a certain minimal projection $q_\alpha \in M$
of our choice (\ref{mark1}) and $u_\beta = u_\gamma v$ for some partial isometry $v \in N$
which links two projections $\{ p_\nu \} \subset N$ of our initial choice.
Then the Hilbert $N$-modules 
$\{ u_\beta (\mu_\beta^{-1/2} N), E(\langle .,. \rangle_M) \}$ and
$\{ u_\gamma (\mu_\gamma^{-1/2} N), E(\langle .,. \rangle_M) \}$
derived from equivalent partial isometries must be unitarily isomorphic, since
the set identity $u_\gamma N \equiv u_\gamma v (v^* N) \equiv u_\beta N$ holds
inside $M$.

\noindent
For every equivalence class we select a representative $u_\beta$ which has to
be a minimal projection $q_\alpha$ of our choice (\ref{mark1}) whenever there
is one contained in the equivalence class under consideration. Since $N$ and
the direct sum of Hilbert $N$-modules
$\,$ w*-$\sum_{\omega \in I} N_{(\omega)} \subseteq M$ are unitarily
isomorphic as Hilbert $N$-modules (where $I$ has the same cardinality as the
index set of the selected sequence of minimal projections $\{ q_\alpha \}
\subset M$), we find that the Hilbert $N$-submodule $u_\beta N$ is unitarily
isomorphic to the direct sum of Hilbert $N$-modules $\,$ w*-$\sum_{\gamma}
u_\gamma N \subseteq M$, where the index runs over all possible indices
$\gamma$ with the property $u_\beta \sim u_\gamma$.

\noindent
Consequently, we can reduce our set of generators (\ref{mark1}) in such a way
that there only remains one element of every equivalence class of partial
isometries, which should be a projection in case the equivalence class
contains one, and which should be a partial isometry connecting
two of the just selected projections otherwise. For further considerations we
will use the same notation as above in (\ref{mark1}) to refer to the reduced
generator set of $M$ generating it as a right Hilbert $N$-module.

\medskip \noindent
After this factorization-like procedure we obtain a possibly smaller generating set of $M$ with
a very special property: for every minimal projection $q_\alpha \in M$ in
this generating set there are at most $([K(E)]-1)$ other minimal projections in
it which are connected to $q_\alpha$ by a partial isometry of this generating
set. Indeed, every minimal projection of this kind has its own minimal
$N$-central carrier projection by construction. Since the W*-factor
$M q_\alpha M$ intersects with at most $[K(E)]$ W*-blocks of the central direct
integral decomposition of $N$ as shown above the statement yields.
To proceed we fix a minimal projection $q_\alpha$ of the actually selected set
that generates $M$ as a Hilbert $N$-module, and we form the sum $p_\alpha$ of all
those minimal projections of $N$ majorizing a minimal projection of $M$ of
our latter choice that is equivalent to $q_\alpha$. This is a finite sum with at most
$[K(E)]$ summands, and the W*-algebra $p_\alpha N p_\alpha$ is commutative.
Therefore, the W*-algebra $p_\alpha M p_\alpha$ is a matrix algebra since
$E$ restricted to it has the commutative finite-dimensional image $p_\alpha N
p_\alpha$, cf.~\cite[Cor.~4.4]{FrKi98}. Finally, the number of modular
generators of our last choice contained in $p_\alpha M p_\alpha$ cannot exceed
$[K(E)]^2$ since the number of minimal projections of $M$ of this choice does
not exceed $[K(E)]$ as shown.

\smallskip \noindent
By transfinite induction we get a partition of the identity of $M$ (and $N$)
as a sum of pairwise orthogonal projections of type $p_\alpha \in N$, and
any generator of the made choice has to be contained in one of the
W*-subalgebras of type $p_\alpha M p_\alpha$ of $M$, where the carrier
projections of these W*-subalgebras are pairwise orthogonal in $N$ by
construction. Consequently, we can form sum-compositions of generators of our
selected minimal set to reduce the number of them further by the following
principle: take one generator per W*-subalgebra of type $p_\alpha M p_\alpha$
and form the appropriate w*-sum of them inside $M$ to get a new generator.
(Equivalently, we can find a subset of singly generated Hilbert $N$-submodules
of our decomposition of $M$ which can be summed up to another singly generated
Hilbert $N$-submodule reducing the number of direct summands that
way enormously.)
In fact, we can form at most $[K(E)]^2$ modular generators of this type since
every W*-subalgebra of type $p_\alpha M p_\alpha$ was shown to contain at most
as many of them of our choice made.

\noindent
Finally, we have constructed a generating set of the Hilbert $N$-module
$\{ M, E(\langle .,. \rangle) \}$ which has at most $[K(E)]^2$ modular
generators, and the theorem is proved.
\end{proof}

\medskip \noindent
{\bf 4. Final remarks. $\,$}
The construction given in the proof works for
properly infinite discrete von Neumann algebras equally well. However, there
is an isomorphism $H \cong \sum_{i \in I} H_{(i)}$ which exists for every
infinite-dimensional Hilbert space $H$ and for ${\rm card}(I) \leq {\rm dim}
(H)$ due to the existence of a chain of pairwise orthogonal projections each
of which is similar to $1_M$ and which sum up to $1_M$. So we can again
reduce our set of generators to a single generator $m \in M$.

\smallskip \noindent
The remaining case to be investigated is the case of finite von Neumann
algebras with diffuse center and some obstruction like $Z(M) \cap Z(N) =
{\bf C} 1_M$, at least valid for parts of the centers. This question is still
unsolved at present.

\medskip \noindent
{\bf Acknowledgements:} I would like to thank C.~Anantharaman-Delaroche,
Y.~Denizeau, F.~Fidaleo, J.-F.~Havet, E.~Kirchberg, V.~Paulsen, S.~Popa and 
J.~Renault for valuable comments and discussions on the subject of this paper.
I am grateful to D.~P.~Blecher and V.~I.~Paulsen for their warm hospitality
and scientific collaboration during a one year stay at the University of
Houston in 1998.


\end{document}